%% file: treekhin3.tex
\documentclass[11pt]{article}
 
\usepackage{amssymb,epsfig} 
 
\input{macros_angl.tex}

\title{A logarithm law for automorphism groups of trees}
 
\author{Sa'ar Hersonsky \and
Fr\'ed\'eric Paulin}

\date{} 
 
\begin{document} 
\maketitle

\begin{abstract} 
\noindent   
Let $\Ga$ be a geometrically finite tree lattice. We prove a
Khintchine-Sullivan type theorem for the Hausdorff measure of the
set of points at infinity of the tree that are well approximated by the
parabolic fixed points of $\Ga$.  Using Bruhat-Tits trees, an
application is given for the Diophantine approximation of formal
Laurent series in the variable $X^{-1}$ over the finite field $\FF_q$
by rational fractions in $X$ over $\FF_q$, satisfying some congruence
properties. $^,\!$
\end{abstract}

\section{Introduction} 
\label{sec:Intro} 
Let $T$ be a locally finite tree, $\partial T$ be its space of ends,
and Aut$(T)$ be its locally compact group of automorphisms. To
simplify the statements in the introduction, we assume that $T$ has no
degree $1$ or $2$ vertex and that $_{{\rm Aut}(T)}\backslash T$ is a
finite graph.  A {\it lattice} of $T$ is a discrete subgroup with a cofinite
Haar measure in Aut$(T)$. It is {\it uniform} if $\Ga\backslash T$ is
a finite graph. In \cite[chap.~4]{BL}, H.~Bass and A.~Lubotzky gave
numerous examples of non uniform lattices. In \cite{Pau2}, the second
author has introduced a special class of discrete subgroups of
Aut$(T)$, where the interesting ones are non uniform lattices, called
geometrically finite (see also Section \ref{sect:definot}). In
particular, they contain all the algebraic examples. More precisely,
let $\hat{K}$ be a non-archimedian local field, and $\underline{G}$ be
a connected semi-simple algebraic group over $\hat{K}$ with
$\hat{K}$-rank $1$.  Let $T$ be the Bruhat-Tits tree of
$(\underline{G},\hat{K})$ (see \cite{BT,Ser}), endowed with its action
of $\underline{G}(\hat{K})$. By the work of J.-P.~Serre,
M.~Raghunathan and A.~Lubotzky \cite{Lub}, every lattice contained in
the image of $\underline{G}(\hat{K})$ in Aut$(T)$ is geometrically
finite in the sense of \cite{Pau2}.

Let $\Ga$ be a fixed geometrically finite lattice in Aut$(T)$. Let
$\Ga_\infty$ be a parabolic subgroup of $\Ga$, $\infty$ be its unique
fixed point in $\partial T$, and $(H_{r})_{r \in\Ga\infty}$ be its
associated family of horoballs (see \cite{Pau2} and also
Section~\ref{sect:definot} for definitions). For every $r$ which is
different from $\infty$ in the orbit $\Ga\infty$, define $D( r)=
d(H_\infty,H_r)$, where $d$ is the distance on $T$.  Let $d_{\partial
T}$ be a usual distance on $\partial T$, $\delta$ be its Hausdorff
dimension and let $\mu$ be the associated Hausdorff measure in
dimension $\delta$ (see Section~\ref{sect:definot}).

\btheo\label{theo:intro-treekhin} %
Let $\psi:\RR_+\ra\RR_+$ be a map such that $\log \psi$ is Lipschitz.
Let $E(\psi)$ be the set of points $\xi$ in $\partial T-\{\infty\}$
such that there exist infinitely many $r$ in $\Ga\infty$ with
$$d_{\partial T}(\xi, r)\leq \psi(D(r))\;e^{-D(r)}\;.$$ Then
$\mu(E(\psi))=0$ (respectively $\mu(^c{}E(\psi))=0$) if and only if
the integral $\int_{1}^{+\infty} \psi(t)^{\delta} dt$ converges
(respectively diverges).  
\etheo

A more general statement is given by Theorem \ref{theo:body-treekhin}.
The main consequence is the following corollary, called the {\it
  logarithm law} for $\Ga$.

Define the height function $h:T\ra [0,+\infty[$ for $\Ga$ by $h(x)= 0$
if $x\in T-\bigcup_{r\in\Ga\infty}H_{r}$, and if $x\in H_{r}$, then
$h(x)= d(x,\partial H_r)$.

\bcoro\label{coro:intro_loglaw} %
For every vertex $x$ in $T$, let $c_\eta$ denote the geodesic ray
starting from $x$ and ending at $\eta\in\partial T$. Then for
$\mu$-almost every $\eta$ in $\partial T$,
$$\limsup_{t\ra+\infty}\frac{h(c_\eta(t))}{\log t}=
\frac{1}{\delta}\;.$$ 
\ecoro

Theorem \ref{theo:intro-treekhin} and this corollary are analogues, in
the setting of trees and in particular for non archimedian rank one
lattices, of the similar well-known results for
\begin{itemize}
 \item 
${\rm PSL}_2(\ZZ)$ acting on the real hyperbolic plane
$\HH^2_\RR$, due to Khintchine \cite{Kh},
\item 
a lattice in the isometry group of a real hyperbolic space
$\HH^n_\RR$, due to D.~Sullivan \cite{Sul},
\item 
a lattice in the isometry group of a symmetric space of non compact
type, due to D.~Kleinbock and G.~Margulis \cite{KM},
\item 
a geometrically finite group of isometries acting on $\HH^n_\RR$ due
 to B.~Stratmann and S.L.~Velani \cite{SV}, and
\item 
more general geometrically finite groups of isometries of complete
simply connected Riemannian manifolds with negative curvature, due to
the authors \cite{HP4}.
\end{itemize}
A particular case of the theorem was announced by J.~Athreya 
\cite{Ath}, when $T$ is the Bruhat-Tits tree of $(\underline{{\rm
    SL}_2}, \FF_q((X^{-1})))$ and $\Ga={\rm SL}_2( \FF_q[X])$.

\medskip %
Denote by $|\cdot|_\infty$ the usual absolute value on the field $\wh
K= \FF_q((X^{-1}))$ of formal Laurent series in $X^{-1}$ over $\FF_q$.
Endowed $\wh K$ with its Haar measure.  Let $Q_0$ be a fixed non zero
element in the ring $A=\FF_q[X]$ of polynomials over $\FF_q$. 
By using \cite{Pau1,Pau2}, we obtain the following corollary,
concerning the Diophantine approximation of elements of $\wh K$ by
rational fractions.

\bcoro\label{coro:intro_diophante} %
Let $\varphi:\RR_+\ra\RR_+$ be a map such that $u\mapsto \log \varphi
(e^u)$ is Lipschitz. If the integral $\int_{1}^{+\infty} \varphi(t)/t
\;dt$ diverges (resp.~converges), then for almost every (resp.~almost
no) $f$ in $\wh K$, there exist infinitely many couples $(P,Q)$ in
$A\times (A-\{0\})$, where $P$ is coprime with $Q$ and $Q$
is divisible by $Q_0$, such that $| f-P/Q |_\infty \leq
\varphi(|Q|_\infty)/|Q|^2_\infty$.  
\ecoro

\medskip %
A more general statement is given by Corollary
\ref{coro:intro_diophante_plus}. The special case of
Corollary~\ref{coro:intro_diophante} when $Q_0$ is the constant
polynomial $1$ was announced by J.~Athreya \cite{Ath}.

\medskip
\noindent {\small {\bf Acknowledgement.} The second author
  acknowledges the support of the University of Georgia for a stay
  during which this paper was written.}

\section{Notations, statements and applications} 
\label{sect:definot}

We refer to \cite{BH,Rob} for basic definitions about CAT$(-1)$
geodesic metric spaces (among which are trees), their horospheres,
their boundaries and their discrete groups of isometries.

Let $T$ be a locally finite tree, $\partial T$ be its space of ends and
$T\cup\partial T$ be its compactification by its ends. Such a tree is
called {\it uniform} if there exists a discrete subgroup $\Ga$ in
Aut$(T)$ such that the quotient graph $\Ga\backslash T$ is finite. The
tree $T$ is endowed with the maximal geodesic distance $d$ making each
edge isometric to the interval $[0,1]$. For every $x,y$ in
$T\cup\partial T $, we denote by $[x,y]$ the geodesic segment, ray or
line between them, with the usual convention concerning the
endpoints. Fix a point $x_0$ in $T$. Let $d_{x_0}$ be the {\it visual
distance} on $\partial T$ seen from $x_0$. It is defined by
$d_{x_0}(\eta,\eta')=e^{-d(x_0,p)}$ where $[x_0,\eta[\, \cap\,
[x_0,\eta'[\;=[x_0,p]$ if $\eta\neq \eta'$. The {\it Buseman function}
$\beta_\xi:X\times X\ra\RR$ of a point $\xi$ in $\partial T$ is defined by
$\beta_\xi(x,y)= d(x,u)-d(y,u)$ for every $u$ in $T$ close enough to
$\xi$. An {\it horoball} centered at $\xi$ is the preimage by
$y\mapsto \beta_\xi(x_0,y)$ of $[t,+\infty[$ for some $t\in \RR$.

Let $\Ga$ be a discrete subgroup of ${\rm Aut}(T)$. We assume that
$\Ga$ is {\it non-elementary}, i.e.~that $\Ga$ preserves no point nor
pair of points in $T\cup\partial T$.  Let $\Lambda \Ga$ be its {\it
  limit set}, i.e.~the smallest non empty closed $\Ga$-invariant
subset in $\partial T$.  Let $C\Lambda \Ga$ be the minimal non empty
subtree of $T$ which is invariant by $\Ga$.

The {\it critical exponent} of $\Ga$ is the unique number $\delta$ in
$[0,+\infty]$ such that the {\it Poincar\'e series} $\sum_{\ga\in\Ga}
e^{-s \;d(x_0,\ga x_0)}$ of $\Ga$ converges for $s>\delta$ and
diverges for $s<\delta$. The group $\Ga$ is called {\it of divergent
  type} if its Poincar\'e series diverges at $s=\delta$.

If $\Ga$ is of divergent type with a finite non zero critical
exponent, then there exists (see \cite{Coo}) a family $(\mu_x)_{x\in
  T}$ of finite measures on $\partial T$, with support
$\Lambda\Gamma$, called the {\it Patterson-Sullivan density}, which is
unique up to a positive scalar factor, such that
\begin{itemize}
\item $\forall\;\ga\in\Ga\;,\;\;\ga_*\mu_x=\mu_{\ga x}$,
\item $\forall\;x,y\in T,\;\forall\;\xi\in\partial
  T\;,\;\;\frac{d\mu_x}{d\mu_y}(\xi)=e^{-\delta \beta_\xi(x,y)}$.
\end{itemize}

An element $\omega$ in $\Lambda \Ga$ is a {\it conical limit point} if
there exist a sequence $(\ga_n)_{n\in\NN}$ in $\Ga$, a geodesic ray
$c$ in $T$ ending at $\omega$ and $A\geq 0$ such that $d(\ga_nx_0,c)\leq
A$ for every $n$ in $\NN$. An element $\omega$ in $\Lambda \Ga$ is
called a {\it bounded parabolic point} if its stabilizer $\Ga_\omega$
in $\Ga$ acts properly discontinuously with compact quotient on
$\partial T-\{\omega\}$. A {\it parabolic subgroup} $\Ga_\infty$ of
$\Ga$ is a maximal infinite locally finite subgroup (see
\cite{BL,Pau2} for other characterizations and the following
properties). In particular, there is one and only one point in
$T\cup\partial T$ fixed by $\Ga_\infty$ and it belongs to $\partial
T$.

A discrete subgroup $\Ga$ of ${\rm Aut}(T)$ is {\it geometrically
finite} if it is non elementary and if every point in $\Lambda \Ga$ is
either a conical limit point or a bounded parabolic point. A
characterization in terms of the quotient graph of groups
$\Ga\backslash\!\backslash T$ (see \cite{BL} for a definition when
$\Ga$ acts without inversion) and the basic properties of
geometrically finite subgroups $\Ga$ of ${\rm Aut}(T)$ are given in
\cite{Pau2}. In particular, the point fixed by a parabolic subgroup of
$\Ga$ is then a bounded parabolic point.

\medskip%
From now on, we fix a geometrically finite subgroup $\Ga$ in Aut$(T)$
and a parabolic subgroup $\Ga_\infty$ of $\Ga$. Let $\infty$ be the
unique fixed point of $\Ga_\infty$ in $\partial T$ and let $\Ga\infty$
be its orbit under $\Ga$.  In \cite{Pau2}, it is proved that there
exists a unique $\Ga$-invariant family of horoballs
$(H_{r})_{r\in\Ga\infty}$, with pairwise disjoint interiors, with
$H_{r}$ centered at $r$ and maximal with respect to inclusion.

Let $d_{\infty}$ be the {\it Hamenst\"adt distance}  on $\partial
T-\{\infty\}$ associated to $H_\infty$ (see \cite[Appendix]{HP1},
\cite{Pau2}). It is defined by
$d_{\infty}(\eta,\eta')=e^{-d_\pm(p,q)}$ where $\{p\}=\partial H_\infty
\cap\,]\infty, \eta[$,
$]\infty,q]=\;]\infty,\eta[\,\cap\,]\infty,\eta'[$ and $d_\pm(p,q)$ is
the signed distance between $p$ and $q$ on $]\infty,\eta[$, if
$\eta\neq\eta'$.  The Hausdorff dimension and the Hausdorff measure
class of $\Lambda\Ga-\{\infty\}$ with respect to the distance
$d_{x_0}$ coincides with the Hausdorff dimension and the Hausdorff
measure class of $\Lambda\Ga-\{\infty\}$ with respect to the
Hamenst\"adt distance $d_\infty$ respectively, as locally on $\partial
T-\{\infty\}$, the distances $d_{x_0}$ and $d_\infty$ differ by a
positive scalar factor.

If $\rho$ is the geodesic ray in $T$ starting from $x_0$ and
converging to $\infty$, then the measures $e^{\delta t}\mu_{\rho(t)}$
on $\partial T-\{\infty\}$ converge, as $t\ra+\infty$, to a measure
$\mu_\infty$ on $\partial T-\{\infty\}$, called the {\it
Patterson-Sullivan measure of $\Ga$ on $\partial T-\{\infty\}$},
supported on $\Lambda\Ga-\{\infty\}$. The measure $\mu_\infty$ has the
same measure class as $\mu_{x_0}$ on $\partial T-\{\infty\}$ (see
\cite{HP4}).

The following theorem summarizes the known results on $\Ga$ that we
will use throughout this paper.

\btheo 
\label{theo:recapprop}
Let $\Ga$ be a geometrically finite subgroup of Aut$(T)$ such that
$C\Lambda\Ga$ is a uniform tree without degree $2$ vertices. Let
$\Ga_\infty$ be a parabolic subgroup of $\Ga$. Then the following
assertions hold.
\begin{enumerate}
\item[(1)] 
  The group $\Ga$ has a finite non zero critical exponent, and 
  is of divergent type. Its action on $\Lambda\Ga$ is ergodic with
  respect to the measure $\mu_{x_0}$. This measure has no atom.
\item [(2)]
There exists a constant $c_1>0$ such that for every $n$ in $\NN$, 
$$
\frac{1}{c_1} e^{n\delta}\leq {\rm Card}\{\ga\in\Ga\;:\; 
d(x_0,\ga x_0)\leq n\}\leq c_1 e^{n\delta}\;.
$$
\item[(3)] 
  The critical exponent $\delta$ equals the Hausdorff
  dimension of $\Lambda\Ga$ with respect to the distance $d_{x_0}$.
\item [(4)] 
  The critical exponent of $\Ga_\infty$ is equal to $\delta/2$.
\item [(5)] 
  There exists a constant $c_2>0$ such that for every $n$ in $\NN$,
$$
\frac{1}{c_2} e^{n\delta/2}\leq {\rm Card}\{\ga\in\Ga_\infty\;:\; 
d(x_0,\ga x_0)\leq n\}\leq c_2 e^{n\delta/2}\;.
$$
\end{enumerate}
\etheo

\dem %
Let $T'=C\Lambda\Ga$. Then $T'$, endowed with the restricted action of
$\Ga$, satisfies the hypotheses of the theorem if and only if $T$,
endowed with the action of $\Ga$, does. If the assertions (1)-(5) are
true for $T'$, then they are also true for $T$. Hence we may assume
that $T=C\Lambda\Ga$. Up to minor changes in the following geometric
series argument when replacing $T$ by its barycentric subdivision,
we may assume that $\Ga$ acts without inversion.

Recall that the quotient graph of groups $\Ga\backslash\!\backslash T$
is the union of a finite graph of finite groups and of finitely many
cuspidal rays of groups (see \cite{BL,Pau2} for the definitions and
\cite{Pau2} for the proof). The covolume of $\Ga$ in Aut$(T)$ is (see
\cite{BL}) a multiple of the sum of the inverses of the cardinals of
the vertex groups in $\Ga\backslash\!\backslash T$.  By an easy
geometric series argument, the fact that $T$ has no vertex of degree
$2$ (here and below, the assumption that there are no degree $2$
vertices inside an horoball of an equivariant family of pairwise
disjoint horoballs centered at all parabolic fixed points is
sufficient) implies in particular that $\Ga$ is a lattice in $T$. (See
for instance \cite{BM}, or simply note that if $\Ga_n$ is the
stabilizer of the $n$-th vertex of a cuspidal ray, then
$[\Ga_{n+1}:\Ga_n]\geq 2$ as $T$ has no vertex of degree $2$, hence
$\sum_n 1/|\Ga_n|\leq \sum_n 1/2^n$ is finite.)

\medskip %
(1) The finiteness (and non vanishing) of the critical exponent
follows for instance from \cite[Prop.~1.7]{BM}. By
\cite[Coro.~6.5]{BM} for instance, the group $\Ga$ is of divergent
type. By \cite[Theorem~1.7]{Rob} and \cite[Corollary~1.8]{Rob} for instance,
the measure $\mu_{x_0}$ is ergodic and has no atom.

(4) See for instance \cite[Prop.~3.1]{BP2}.

(2) As every vertex of $T$ has degree at least $3$, and as there
exists at least one parabolic subgroup, the subgroup of $\ZZ$
generated by the hyperbolic translation lengths of the elements of
$\Ga$ equals to $\ZZ$ (see \cite{BP2}). By (1), (4) and
\cite[Th\'eor\`eme~1.11]{Rob}, the Bowen-Margulis measure of $\Ga$ (see for
instance \cite{Rob} for a definition of this measure) is finite. Hence
by corollaire 2 following th\'eor\`eme 4.1.1 in \cite{Rob}, the result
holds.

(3) This follows for instance from \cite[Coro.~6.5]{BM}.

(5) This follows for instance from \cite[Prop.~3.1]{BP2}.
\cqfd

\bigskip
Before providing a more general statement than Theorem
\ref{theo:intro-treekhin}, we will introduce the class of functions we
will work with.

A map $\psi:\RR_+\ra\RR_+$ is called {\it slowly varying} (see
\cite{Sul}) if it is measurable and if there exist constants $B>0$ and
$A\geq 1$ such that for every $x,y$ in $\RR_+$, if $|x-y|\leq B$, then
$\psi(y)\leq A\,\psi(x)$.  Recall (see for instance
\cite[Sec.~5]{HP4}) that this implies that $\psi$ is locally bounded,
hence it is locally integrable; also, if $\log \psi$ is Lipschitz,
then $\psi$ is slowly varying. 

\btheo\label{theo:body-treekhin} 
Let $\Ga$ be a geometrically finite subgroup of Aut$(T)$ such that
$C\Lambda\Ga$ is a uniform tree without degree $2$ vertices. Let
$\Ga_\infty$ be a parabolic subgroup whose fixed point will be denoted
by $\infty$.  Let $\psi:\RR_+\ra\RR_+$ be a slowly varying map.  Let
$E(\psi)$ be the set of points $\xi$ in $\partial T-\{\infty\}$ such
that there exist infinitely many $r$ in $\Ga\infty$ with
$d_{\infty}(\xi, r)\leq \psi(D(r))e^{-D(r)}$. Then
$\mu_\infty(E(\psi))=0$ (respectively $\mu_\infty(^c{}E(\psi))=0$) if
and only if the integral $\int_{1}^{+\infty} \psi(t)^{\delta} dt$
converges (respectively diverges).  
\etheo

The proof of this theorem will be given in Section
\ref{sect:proof}. Note that for some of the Bruhat-Tits trees (the
ones that are biregular of type $(2,q)$, see \cite{BT}), we have to
remove the vertices of degree $2$ in order to verify the hypothesis of
the above theorem. In the new tree, the critical exponent is
multiplied by $2$, the complexities $D(r)$'s are divided by $2$ and
the distance $d_\infty$ is replaced by its square root. Hence the
conclusion of Theorem~\ref{theo:body-treekhin} is also true for these
Bruhat-Tits trees.

Let us show now that Theorem \ref{theo:intro-treekhin} and Corollary
\ref{coro:intro_loglaw} in the introduction follows from
Theorem~\ref{theo:body-treekhin}.

\medskip
\noindent{\bf Proof of Theorem \ref{theo:intro-treekhin}.} %
Since $\Ga$ in the statement of Theorem \ref{theo:intro-treekhin} is
assumed in addition to be a lattice, we have $\Lambda\Ga=\partial T$.
Hence, since $T$ has no terminal vertex, $C\Lambda\Ga=T$. A tree whose
quotient by its automorphism group is a finite graph, and whose
automorphism group contains a (non-uniform) lattice is unimodular,
hence uniform (see \cite[para.~1.2, 0.5 (6)]{BL}).  As we noticed above, a
map whose logarithm is Lipschitz is slowly varying. Hence, if the
hypotheses of Theorem
\ref{theo:intro-treekhin} are satisfied, then so are the hypotheses of
Theorem~\ref{theo:body-treekhin}. The critical exponent of $\Ga$
coincides with the Hausdorff dimension of $\Lambda\Ga$ with respect to
the $d_{x_0}$ metric.  The distances $d_\infty$ and $d_{x_0}$ are
locally equivalent, and the measure classes of the measure
$\mu_\infty$ and of the Hausdorff measure of $d_{x_0}$ coincide. Hence
Theorem \ref{theo:intro-treekhin} does follow from
Theorem~\ref{theo:body-treekhin}.  
\cqfd

\bigskip
\noindent{\bf Proof of Corollary \ref{coro:intro_loglaw}.} %
The proof follows from Theorem 1.1 along the same line as in the other
cases we recalled in the introduction (see for instance
\cite[Sect.~6]{HP4}).  We use the maps $\psi(t)=t^{-\kappa}$, and the
fact that for a geodesic line starting from $\infty$, ending at $\xi$,
and entering in $H_r$ for some $r\in\Ga\infty-\{\infty\}$ with highest
penetration point $\xi_r$, we have $-\log d_\infty(\xi,r)=
D(r)+d(\xi_r,\partial H_r)$.  
\cqfd

\bigskip %
We now give some applications of our results.  We refer for instance
to \cite{Las,Sch} for nice surveys of the Diophantine approximation
properties of elements in $\wh K= \FF_q((X^{-1}))$ by elements in
$K=\FF_q(X)$.  Recall the definition of the absolute value of $f\in\wh
K-\{0\}$.  Let $f=\sum_{i=n}^\infty a_i X^{-i}$ where $a_n\neq 0$.
Then we define $\nu(f)=n\in \ZZ$ and $|f|_\infty=q^{-\nu(f)}$.  Endow
the locally compact additive group $\wh K$ with its (unique up to a
constant factor) Haar measure. Let $A=\FF_q[X]$.

\bcoro\label{coro:intro_diophante_plus} %
Let $\Ga$ be a finite index subgroup of ${\rm SL}_2(A)$.  Let
$\varphi:\RR_+\ra\RR_+$ be a map with $u\mapsto \varphi (e^u)$ slowly
varying. If the integral $\int_{1}^{+\infty} \varphi(t)/t \;dt$
diverges (resp.~converges), then for almost every (resp.~almost no)
$f$ in $\wh K$, there exist infinitely many couples $(P,Q)$ in
$A\times (A-\{0\}) $, where $P$ is coprime with $Q$ and
$P/Q\in\Ga\infty$, such that $| f-P/Q |_\infty \leq
\varphi(|Q|_\infty)/|Q|^2_\infty$.  
\ecoro

\dem In \cite{Pau1}, a geometric interpretation of the above
Diophantine approximation is given in terms of the Bruhat-Tits tree
$\TT_q$ of $({\underline{\rm SL}}_2, \wh K)$ (see also
\cite{Pau2,BP1}).  Identify as usual $\partial \TT_q$ and $\wh
K\cup\{\infty\}$, and let $x_0$ be the standard base point in $\TT_q$
(see \cite{Ser}). Note that the Hausdorff dimension of $d_{x_0}$ is
$\log q$, as $\TT_q$ is a regular tree of degree $q+1$ (see
\cite{Ser}).  Recall that the kernel of the action of ${\rm SL}_2(\wh
K)$ on $\TT_q$ is finite, hence denoting by the same letter a lattice
in ${\rm SL}_2(\wh K)$ and its image in ${\rm Aut}(\TT_q)$ causes no
problem.

Now $\Ga$ is a lattice of $\TT_q$, hence $C\Lambda\Ga=\TT_q$ is a
uniform tree without degree $2$ vertices. Let $\Ga_\infty$ be the
stabilizer in $\Ga$ of the point $\infty$ in $\partial \TT_q=\wh
K\cup\{\infty\}$, which is a parabolic subgroup of $\Ga$. Since $\Ga$
is contained in ${\rm SL}_2(A)$, the subset $\Ga\infty-\{\infty\}$ is
contained in $K$. By \cite{Pau1,Pau2}, the following assertions hold.
\begin{itemize}
\item[$\bullet$] 
  There exists a constant $c>0$ such that the Hamenst\"adt distance 
  on $\partial \TT_q-\{\infty\}=\wh K$ satisfies
  $d_\infty(f,f')=c|f-f'|_\infty^{\frac{1}{\log q}}$, for every $f,f'$
  in $\wh K$ (see \cite[Coro.~5.2]{Pau1}).
\item[$\bullet$] 
  There exists a constant $c'>0$ such that for every
  $P/Q$ in $\Ga\infty-\{\infty\}$, where $P$ and $Q$ are coprime
  polynomials such that $\nu(P)\geq\nu(Q)$, we have $D(P/Q)=-2\nu(Q)
  +c'=(2\log |Q|_\infty)/\log q +c'$ (see \cite[Coro.~6.1]{Pau1}).
\end{itemize}
The constants are due to the fact that for $\ga\in \Ga$, the horoball
centered at a point $\ga\infty$ in the family associated to $\Ga$
contains, but is not necessarily equal to, the horoball centered at a
point $\ga\infty$ in the family associated to ${\rm SL}_2(A)$.  Hence,
for every $c''>0$, we have $d_{\infty}(f, P/Q)\leq c''e^{-D(P/Q)}$ if
and only if $|f-P/Q|_\infty\leq (\frac{c''e^{-c'}}{c})^{\log
q}/|Q|^2_\infty$. Define $\psi(t)= ce^{c'}\;\varphi(q^{(t-c')/2})
^\frac{1}{\log q}$, for every $t>0$. So that $\psi$ is slowly varying
by the assumption on $\varphi$. And $d_{\infty}(f, P/Q)\leq
\psi(D(P/Q))e^{-D(P/Q)}$ if and only if $|f-P/Q|_\infty\leq
\varphi(|Q|_\infty)/|Q|^2_\infty$.

By an easy change of variables, the integral $\int_{1}^{+\infty}
\varphi(t)/t \;dt$ diverges if and only if $\int_{1}^{+\infty}
\psi(t)^{\log q}\,dt$ diverges.  By invariance, the Haar measure of
$\wh K$ is equal to the Hausdorff measure of $d_\infty$, up to a
constant positive factor. In particular, their measure classes are the
same. The result follows.  
\cqfd

\medskip
\noindent{\bf Proof of Corollary \ref{coro:intro_diophante}.}  %
Apply Corollary \ref{coro:intro_diophante_plus} to the congruence
subgroup $\Ga$ of ${\rm SL}_2(A)$, consisting of the elements
$\left(\begin{array}{cc}a & b \\ c & d \end{array}\right)$ such that
$c$ is divisible by $Q_0$. The subgroup $\Ga$ has finite index in
${\rm SL}_2(A)$, as it is the preimage of the upper triangular
subgroup by the reduction modulo $Q_0$ map:
$$
\;\;\;\;{\rm Ker}({\rm SL}_2( \FF_q[X])\ra {\rm SL}_2(
\FF_q[X]/Q_0\,\FF_q[X]))\;,
$$ with the convention that ${\rm SL}_2(R)$ is the trivial group if
$R$ is the trivial one-point ring. Note that $\Ga\infty-\{\infty\}$ is
the set of fractions $P/Q$ where $P$ is coprime with $Q$ and $Q$ is
divisible by $Q_0$.  
\cqfd

\medskip Similar results to Corollary \ref{coro:intro_diophante} can
be obtained for instance by varying the congruence subgroup of ${\rm
  SL}_2( \FF_q[X])$.

\section{Proof of the main result}
\label{sect:proof}
Throughout this section, we will assume that the hypotheses of Theorem
\ref{theo:body-treekhin} are satisfied. As $E(\psi)$ is contained in
$\Lambda\Ga$, and as $\mu_\infty$ has support in $\Lambda\Ga$, we may
assume, up to replacing $T$ by $C\Lambda\Ga$, that $T=C\Lambda\Ga$.

We first claim that for every constant $\eta>0$, in order to prove
Theorem \ref{theo:body-treekhin}, it is sufficient to prove the
assertion if $\psi(t)\leq \eta$ for every $t$ in $\RR_+$.  Indeed,
since $\infty$ is a bounded parabolic point and $\Ga$ acts discretely
on $T$, there exist, modulo the action of the parabolic subgroup
$\Ga_\infty$, only finitely many elements $r$ in
$\Ga\infty-\{\infty\}$ with $D(r)$ less than some constant.  The above
reduction is then justified by the same arguments as in the proof of
the similar reduction in \cite[Lem.~5.2]{HP4}. In particular, from now
on, we will assume that $\psi\leq 1$.

\medskip
Let $\Delta$ be a compact subset of $\partial T-\{\infty\}$ whose
images under $\Ga_\infty$ cover $\partial T-\{\infty\}$, and whose
interior contains a representative of every element of
$\Ga\infty-\{\infty\}$ modulo $\Ga_\infty$. The existence of such a
subset $\Delta$ follows from the fact that $\infty$ is a bounded
parabolic point.  The following proposition gives counting estimates
on the number of points in $(\Ga_\infty)\cap\Delta$.

\bprop \label{prop:depthgrowth} 
There exists an integer $N\geq 1$ and a constant $c_3>0$ such that for
 every $n$ in $\NN$, if $\N(n)$ denotes the number of elements $r$ in
 $(\Ga\infty)\cap\Delta$ such that $n\leq D(r) < n +N$, then
$$\frac{1}{c_3}e^{\delta n}\leq \N(n)\leq c_3e^{\delta  n}\;.$$ 
\eprop

\dem %
Let $\N'(n)= {\rm Card} \{[\ga]\in \Ga_\infty\backslash
(\Ga-\Ga_\infty)/\Ga_\infty\;:\; d(H_\infty,\ga H_\infty)\leq n\}$.
We first claim that there exists a constant $c'_1>0$ such that
$e^{\delta n}/c'_1\leq \N'(n)\leq c'_1e^{\delta n}$. The proof of this
claim is the same (with simplifications due to the tree structure) as
the proof of \cite[Theo.~3.4]{HP4}. The crucial ingredients we used in
that proof were that $\infty$ is a bounded parabolic point, that the
orbit $\Ga x_0$ satisfies some growth property (which is given here by
Theorem \ref{theo:recapprop} (2)) and that the critical exponent of
$\Ga_\infty$ is strictly less than $\delta$ (this is given here by
Theorem \ref{theo:recapprop} (4)).

By the properties of $\Delta$, there exists a constant $c>0$ such that
$$\N'(n)\leq{\rm Card} \{r\in (\Ga\infty)\cap\Delta\;:\;D(r) \leq n\}
\leq c\;\N'(n)\;.$$ Hence there exists $c''_1>0$ such that $e^{\delta
  n}/c''_1\leq {\rm Card} \{r\in (\Ga\infty)\cap\Delta\;:\;D(r) \leq
n\}\leq c''_1e^{\delta n}\;.$ Finally, Proposition
\ref{prop:depthgrowth} follows from this, and from
\cite[Lem.~3.3]{HP4}, as in the proof of \cite[Theo.~3.4]{HP4}.  
\cqfd

\bigskip %
If $H$ is an horoball and $t\geq 0$, define $H(t)$ to be the horoball
contained in $H$, whose boundary is at distance $t$ from the boundary
of $H$. To simplify the notation, for every $r$ in
$\Ga\infty-\{\infty\}$, let $H_{r,\psi}=H_r(-\log \psi\circ D(r))$.

For every subset $E$ of $T$, let $\O_\infty E$ be the set of endpoints
of the geodesic lines starting from $\infty$ and passing through $E$.
Let $A_n$ be the subset of $\partial T-\{\infty\}$, which is the union
of the  $\O_\infty H_{r,\psi}$'s where $r$ ranges over the elements
in $(\Ga\infty)\cap\Delta$ with $Nn\leq D(r)< (n+1)N$ and $N$ is as in
Proposition \ref{prop:depthgrowth}.  

\medskip
The proof of Theorem \ref{theo:body-treekhin} is based on the
following three propositions.  The first one is an adaptation for
trees of the fluctuating density property due to Sullivan \cite{Sul}
in the case of finite volume, constant curvature manifolds, to
Stratmann-Velani \cite{SV} in the geometrically finite, constant
curvature manifold case, and to \cite{HP4} for more general
geometrically finite manifolds with variable negative curvature.  The
proof of the following proposition greatly simplifies in the case of
homogeneous trees, but the general case requires some care.

\bprop \label{prop:fluctuat_density} %
There exists a constant $c>0$ such that for every $r$ in
$\Ga\infty-\{\infty\}$ and every $t\geq 0$, one has
$$\frac{1}{c}e^{-\delta (D(r) + t)}\leq 
\mu_{\infty}(\O_\infty H_{r}(t))\leq c e^{-\delta (D(r) + t)}\;.$$ 
\eprop

\dem %
It is the same as the proof of \cite[Theo.~4.1]{HP4}, except for the
second step of the proof of \cite[Prop.~4.3]{HP4}. To check the
hypotheses of \cite[Theo.~4.1]{HP4}, we use the following claims: the
point $\infty$ is a bounded parabolic point, $\Ga$ is of divergent
type (see Theorem \ref{theo:recapprop} (1)) and for every $n$ in
$\NN$,
$$\frac{1}{c_2} e^{n\delta/2}\leq {\rm Card}\{\ga\in\Ga_\infty\;:\;
d(x_0,\ga x_0)\leq n\}\leq c_2 e^{n\delta/2},$$ (see Theorem
\ref{theo:recapprop} (5)), so that the constant $\delta_0$ that appears
in the statement of \cite[Theo.~4.1]{HP4} is now $\delta/2$.
Furthermore, Proposition 4.2 in \cite{HP4}, which is needed in the
proof of Theorem 4.1 in \cite{HP4}, is a bit simpler in the setting of
trees.

Step 2 in the proof of \cite[Prop.~4.3]{HP4} (which is also needed for
the proof of \cite[Theo.~4.1]{HP4}) needs to be adjusted.  The
arguments that used the hypothesis of pinched curvature are no longer
valid. But, with the notation of Step 2 in the proof of
\cite[Prop.~4.3]{HP4}, the convergence of
$(X_i,\ast_i,d_i,G_i)_{i\in\NN}$ follows from the fact that these
spaces are uniformly proper (see \cite{Gro}), as $T$ is a uniform
tree. The limit is again a locally finite tree, and the rest of the
proof is unchanged.  \cqfd

\bprop\label{prop:khinchine_deux}  
The sum $\displaystyle \sum_{n=0}^{\infty} \mu_\infty(A_n)$ 
diverges if and only if the integral $\displaystyle\int_{1}^{\infty} 
\psi{\,}^{\delta} $ diverges. 
\eprop 
 
\dem The proof of this proposition is completely similar to the proof
of \cite[Prop.~5.3]{HP4}. Lemma 5.4 used in that proof is first proved
for the case of trees.  Proposition \ref{prop:fluctuat_density} above
replaces \cite[Theo.~4.1]{HP4}.  
\cqfd

\bprop\label{prop:khinchine_un} 
There exists a constant $c>0$ such that if $n,m$ are distinct 
 integers, then 
$$\mu_\infty(A_n\cap A_m)\leq c\mu_\infty(A_n)\mu_\infty(A_m).$$ 
\eprop  

\dem The proof of this proposition is completely similar to the proof
of \cite[Prop.~5.5]{HP4}.  Lemma 5.6 used in that proof is even first
proved for the case of trees.  Proposition \ref{prop:fluctuat_density}
above replaces \cite[Theo.~4.1]{HP4}. \cqfd

\bigskip %
We now end the proof of Theorem \ref{theo:body-treekhin}.  Recall that
in any tree $T'$, given an horoball $H$ with point at infinity $\xi$,
and distinct points $\eta,\eta'\in\partial T'$ different from $\xi$,
if $x$ is the intersection with $\partial H$ of the geodesic line
$]\eta, \xi[$, then $\eta'$ belongs to $\O_\eta H$ is and only if the
geodesic line $]\eta, \eta'[$ goes through $x$.

In particular, as the distance between $H_\infty$ and $H_{r,\psi}$ is
equal to $D(r)-\log \psi\circ D(r)$, a point $\xi$ in $\partial
T-\{\infty\}$ belongs to $\O_\infty H_{r,\psi}$ if and only if
$d_\infty(\xi,r)\leq \psi\circ D(r)\;e^{-D(r)}$.

Define $A_\infty=\bigcap_{n\in\NN}\bigcup_{k\geq n} A_k$, which is the
set of points in $\partial T-\{\infty\}$ belonging to infinitely many
$A_n$'s.  With the notation in the statement of Theorem
\ref{theo:body-treekhin}, as the images of $\Delta$ by $\Ga_\infty$
cover $\partial T-\{\infty\}$, up to enlarging $\Delta$, we then have
$E(\psi)= \Ga_\infty A_\infty$. As $\Ga_\infty$ is countable, we have
$\mu_\infty(E(\psi))>0$ if and only if $ \mu_\infty(A_\infty)>0$.

\medskip The following result is well known (see for instance
\cite{Spr}).
 
\bprop\label{prop:borel_cantelli}  
Let $(Y,\nu)$ be a measurable space with a finite measure. 
Let $(B_n)_{n\in\NN}$ be a 
sequence of measurable subsets of $Y$ such that there exists a 
constant $c>0$ with $\nu(B_n\cap B_m)\leq c\nu(B_n)\nu(B_m)$ for every 
distinct integers $n,m$.  Let $ 
B_\infty=\bigcap_{n\in\NN}\bigcup_{k\geq n} B_k$. Then 
$\nu(B_\infty)>0$ if and only if $\displaystyle \sum_{n=0}^{\infty} 
\nu(B_n)$ diverges. \hfill\cqfd  
\eprop

Let $Y$ be a big enough compact subset of $\partial T-\{\infty\}$
which contains $\bigcup_{n\in\NN} A_n$. Let $\nu$ be the restriction
of $\mu_\infty$ to $Y$, and let $B_n=A_n$.  We now apply Proposition
\ref{prop:khinchine_un}, Proposition \ref{prop:khinchine_deux}, and
the result above, to obtain that $\mu_\infty(E(\psi))>0$ if and only
if $\int_{1}^{\infty} \psi^{\delta}$ diverges. This gives the first
conclusion of Theorem \ref{theo:body-treekhin}.

Conversely, assume that $\int_{1}^{\infty} \psi^{\delta}$ diverges.
Let us prove that $\mu_\infty(^cE(\psi))=0$.

Let $g:[0,+\infty[\,\ra\,]0,+\infty[$ be a map decreasing to $0$ such
that $\int_{1}^{\infty} (g\psi)^{\delta}$ diverges. Let $E'(\psi)$ be
the set of $\xi$ in $\partial T-\{\infty\}$ such that there exist
$c>0$ and infinitely many $r$ in $\Ga\infty-\{\infty\}$ with
$d_\infty(\xi,r)\leq c \psi( D(r))e^{-D(r)}$.  Since $E(g\psi)\subset
E'(g\psi)$, the first conclusion of Theorem \ref{theo:body-treekhin}
implies that $\mu_\infty(E'(g\psi))>0$. It is clear that $E'(g\psi)$
is a measurable subset of $\partial T$ which is invariant under
$\Gamma$.

By Theorem \ref{theo:recapprop} (1), the action of $\Ga$ on $\partial
 T$ for the measure $\mu_{x_0}$ is ergodic and without atom.  By
 ergodicity, $\mu_\infty(^cE'(g\psi))=0$.  But $E'(g\psi)\subset
 E(\psi)$ since $g$ is decreasing to $0$. Hence
 $\mu_\infty(^cE(\psi))=0$.  This ends the proof of Theorem
 \ref{theo:body-treekhin}.\cqfd

\noindent {\small 
\begin{tabular}{l}  
University of Georgia\\
Department of Mathematics\\
Athens, GA 30602, USA\\
{\it e-mail: saarh@math.uga.edu} 
\end{tabular} \hfill and \hfill 
\begin{tabular}{l}  
Ben-Gurion University\\
Department of Mathematics\\
Beer-Sheva, Israel\\
{\it e-mail: saarh@math.bgu.ac.il}
\end{tabular}
\\ 
 \mbox{} 
\\ 
 \mbox{} 
\\ 
\begin{tabular}{l}  
D\'epartement de Math\'ematique et Applications, UMR 8553 CNRS\\ 
\'Ecole Normale Sup\'erieure, 45 rue d'Ulm\\ 
75230 PARIS Cedex 05, FRANCE\\ 
{\it e-mail: Frederic.Paulin@ens.fr} 
\end{tabular} 
}

\end{document}

%% file: macros_angl.tex
\newtheorem{defi}{Definition}[section]
\newtheorem{prop}[defi]{Proposition}
\newtheorem{theo}[defi]{Theorem}
\newtheorem{conj}[defi]{Conjecture}
\newtheorem{lemm}[defi]{Lemma}
\newtheorem{coro}[defi]{Corollary}
\newtheorem{rema}[defi]{Remark}
\newtheorem{exem}[defi]{Example}
\newtheorem{exems}[defi]{Examples}
\newtheorem{prob}[defi]{Problem}

\newcommand{\bdefi}{\begin{defi}}
\newcommand{\edefi}{\end{defi}}
\newcommand{\bprop}{\begin{prop}}
\newcommand{\eprop}{\end{prop}}
\newcommand{\btheo}{\begin{theo}}
\newcommand{\etheo}{\end{theo}}
\newcommand{\blemm}{\begin{lemm}}
\newcommand{\brema}{\begin{rema}}
\newcommand{\erema}{\end{rema}}
\newcommand{\bexer}{\begin{exem}}
\newcommand{\eexer}{\end{exem}}
\newcommand{\bexems}{\begin{exems}}
\newcommand{\eexems}{\end{exems}}
\newcommand{\bconj}{\begin{conj}}
\newcommand{\econj}{\end{conj}}
\newcommand{\elemm}{\end{lemm}}
\newcommand{\bcoro}{\begin{coro}}
\newcommand{\ecoro}{\end{coro}}
\newcommand{\bpro}{\begin{prob}}
\newcommand{\epro}{\end{prob}}
\newcommand{\dem}{\noindent{\bf Proof. }}

\newcommand{\M}{{\cal M}}
\newcommand{\N}{{\cal N}}

\renewcommand{\O}{{\cal O}}

\newcommand{\maths}[1]{{\mathbb #1}}  

\newcommand{\RR}{\maths{R}}
\newcommand{\NN}{\maths{N}}

\newcommand{\HH}{\maths{H}}
\newcommand{\FF}{\maths{F}}

\newcommand{\ZZ}{\maths{Z}}

\newcommand{\TT}{\maths{T}}

\newcommand{\ra}{\rightarrow}

\newcommand{\wh}[1]{{\widehat{#1}}}

\newcommand{\ga}{\gamma}
\newcommand{\Ga}{\Gamma}

\newcommand{\cqfd}{\hfill$\Box$}

\newcounter{fig}

